\documentclass[11pt,a4paper]{amsart}
\usepackage[a4paper,hscale=0.65,vscale=0.65,centering]{geometry}

\newcommand{\E}{\mathbb{E}}
\renewcommand{\P}{\mathbb{P}}
\newcommand{\erre}{\mathbb{R}}

\newcommand{\dbr}[1]{\langle\!\langle#1\rangle\!\rangle}
\newcommand{\ds}{\displaystyle}
\newcommand{\ip}[2]{\left\langle#1,#2\right\rangle}
\newcommand{\tqv}[1]{[\![#1]\!]}
\newcommand{\tr}{\mathop{\mathrm{Tr}}\nolimits}

\newtheorem{proposition}{Proposition}[section]
\newtheorem{theorem}[proposition]{Theorem}
\newtheorem{lemma}[proposition]{Lemma}

\newtheorem{definition}[proposition]{Definition}
\theoremstyle{remark}
\newtheorem{remark}[proposition]{Remark}

\numberwithin{equation}{section}

\begin{document}

\title[Stochastic porous media equations with jumps]{Strong solutions
  for stochastic porous media equations with jumps}

\author{Viorel Barbu}
\address{University Al.~I.~Cuza, 8 Blvd. Carol I, Ia\c{s}i 700506, Romania.}
\email{vb41@uaic.ro}
\author{Carlo Marinelli}
\address{Institut f\"ur Angewandte Mathematik, Universit\"at Bonn,
  Wegelerstr. 6, D-53115 Bonn, Germany.}
\urladdr{http://www.uni-bonn.de/$\sim$cm788}

\begin{abstract}
  We prove global well-posedness in the strong sense for stochastic
  generalized porous media equations driven by a square
  integrable martingale with stationary independent increments.
\end{abstract}



\thanks{A large part of this work was written while the second-named author
was visiting the Departments of Mathematics and Statistics of Purdue
University supported by a MOIF grant of the EU}

\maketitle


\section{Introduction}
The purpose of this note is to establish well-posedness in the
strong sense for a class of nonlinear stochastic PDEs driven by L\'evy
noise. More precisely, we shall consider the following stochastic
porous media equation
\begin{equation}
\label{eq:spm}
dX(t) - \Delta\beta(X(t))\,dt = B(X(t-))\,dM(t),
\end{equation}
where $M$ is a square integrable martingale with stationary
independent increments taking values in a Hilbert space $K$, and the
diffusion coefficient $B$ satisfies a Lipschitz assumption.  Full
details on the data of the problem are given below. The present paper
is a continuation of and should be read together with
\cite{BDPR-porous}, where (\ref{eq:spm}) is studied assuming that $M$
is a Wiener process.

Let us just briefly mention that the deterministic counterpart of
(\ref{eq:spm}), i.e. with $B \equiv 0$, has been extensively studied
both for its physical importance and as a model nonlinear PDE (see
e.g. \cite{DaskaKen-libro,VazPME} for systematic treatments). Porous
media equations perturbed by Wiener noise have also been intensively
investigated in the past few years (see references in
\cite{BDPR-porous}).

It should be said that without polynomial growth assumptions on
$\beta$, one cannot apply variational methods (see
e.g. \cite{Gyo-semimg,MP}), and since the drift contains no linear
term, the semigroup approach does not apply either (see
e.g. \cite{PZ-libro}). We have also been unable to find in the
literature results on well-posedness in the strong or mild sense for
nonlinear SPDEs with discontinuous noise that do not fall into any of
the two mentioned settings.

Here we show that (\ref{eq:spm}) admits a unique solution
which depends continuously on the initial datum, thus extending the
results of \cite{BDPR-porous} to the case of a jump noise. Moreover,
we prove that the solution of (\ref{eq:spm}) lives in a ``better''
space than the one used in \cite{BDPR-porous}, and we introduce a
concept of generalized solution that allows to remove some
restrictions on the coefficient $B$ used in \cite{BDPR-porous}.

Let us conclude this introductory section with a few words about
notation. Given two separable Hilbert spaces $H$, $K$ we shall denote
the space of Hilbert-Schmidt operators from $H$ to $K$ by
$\mathcal{L}_2(H,K)$, and the space of trace-class operators on $H$ by
$\mathcal{L}_1(H)$. Moreover, $\mathcal{L}_1^+$ stands for the subset
of $\mathcal{L}_1$ of positive operators. Given a self-adjoint
operator $Q\in\mathcal{L}_1^+(H)$, we denote by $\mathcal{L}_2^Q$ the
set of all (possibly unbounded) operators $B:Q^{1/2}H\to K$ such that
$BQ^{1/2}\in\mathcal{L}_2(H,K)$. The norm in $\mathcal{L}_2^Q$ will be
denoted by $|\cdot|_Q$. We shall denote the space of weakly continuous
functions defined on the interval $I \subseteq \erre$ and taking
values in a Banach space $X$ by $C^w(I,X)$.
Throughout the paper, $\Xi$ will be an open bounded subset of
$\erre^d$, $d>1$, with smooth boundary $\partial\Xi$, and 
\[
Q_t := (0,t) \times \Xi,
\qquad
\partial Q_t := (0,t) \times \partial\Xi,
\qquad t>0.
\]
We denote by $H^{-1}:=H^{-1}(\Xi)$ the dual of the Sobolev space
$H_0^1(\Xi)$, endowed with the scalar product
$\ip{f}{g}_{-1}=\ip{(-\Delta)^{-1}f}{g}$, where $\ip{\cdot}{\cdot}$ is
the duality pairing between $H_0^1(\Xi)$ and $H^{-1}(\Xi)$, and
$\Delta$ stands for the Laplacian with Dirichlet homogeneous boundary
conditions. The inner product in $L^2(\Xi)$ will be denoted by
$\langle\cdot,\cdot\rangle_2$. Whenever no misunderstanding can arise,
we shall write all functional spaces without explicitly indicating the
domain $\Xi$, e.g. $H^{-1} = H^{-1}(\Xi)$, etc.


\section{Main result}
Let us first recall a few facts from the theory of stochastic
integration in Hilbert spaces. For all unexplained notation we refer
to \cite{Met}. Denoting the space of locally square
integrable martingales on $K$ by $\mathcal{M}_{loc}^2(K)$, let $M \in
\mathcal{M}_{loc}^2(K)$ and $T>0$ a fixed (deterministic)
time. Appealing to \cite[{\S}22]{Met}, one can construct stochastic
integrals of the type
\begin{equation}
  \label{eq:sint}
G\cdot M(t) := \int_{(0,t]} G(s)\,dM(s), \qquad t \in [0,T],  
\end{equation}
for a class of operator valued predictable processes $G$.
In particular, according to the results of \cite{Met,MP-Z},
there exists a predictable $\mathcal{L}_1^+(K)$-valued process
$(Q_M(t))_{t\in[0,T]}$ such that
$$
\dbr{M}(t) = \int_0^t Q_M(s)\,d\langle M\rangle(s),
\qquad t\in [0,T],
$$
and the stochastic integral (\ref{eq:sint}) is well defined for all
predictable processes $G:[0,T]\times\Omega \to \mathcal{L}_2^{Q_M}(K,H)$
such that
\begin{equation}
  \label{eq:G}
\E\int_0^T |G(s)Q_M(s)^{1/2}|^2_{\mathcal{L}_2(K,H)}\,d\langle M\rangle(s)
< \infty.
\end{equation}
The set of all such processes will be denoted by $\mathcal{G}(H)$.
Here and everywhere in the following we shall write, with a slight
abuse of notation, $\int_0^t$ instead of $\int_{(0,t]}$.
\begin{remark}     \label{rmk:pissi}
  If $M$ has stationary independent increments (i.e. $M$ is also a
  L\'evy process), then $Q_M=(\tr Q)^{-1}Q$, where $Q$ is the
  covariance operator of $M$, and $\langle M\rangle(t)=t\tr Q$,
  $\dbr{M}(t)=t(\tr Q)^{-1}Q$ for all $t\in[0,T]$. Moreover, $Q$ is a
  \emph{deterministic} operator.
\end{remark}
We shall study existence, uniqueness and regular dependence on the initial
datum for the following stochastic Cauchy problem:
\begin{equation}
\label{eq:main}
\begin{cases}
dX(t)-\Delta\beta(X(t))\,dt = B(X(t-))\,dM(t) &\text{in $Q_T$}\\
\beta(X(t))=0 &\text{in $\partial Q_T$}\\
X(0)=x &\text{in $\Xi$}
\end{cases}
\end{equation}
under a set of assumptions on $\beta$, $B$ and $M$ precised below.
In particular, we shall assume that the diffusion coefficient is of
the form
\begin{equation}     \label{eq:asso}
B:H^{-1} \to \mathcal{L}_2^Q(K,D((-\Delta)^\gamma),
\qquad \gamma>d/2,
\end{equation}
and satisfies
\begin{equation}      \label{eq:assb}
|B(x)|_Q^2 \leq k(1+|x|^2),
\quad
|B(x)-B(y)|_Q^2 \leq k|x-y|^2,
\end{equation}
for some constant $k>0$. Assumption (\ref{eq:asso}) will be relaxed in
the last section.

The (multivalued) function $\beta:\erre \to 2^\erre$ is assumed to be
a maximal monotone graph in $\erre \times \erre$ such that $0 \in
\beta(0)$, $D(\beta)=R(\beta)=\erre$, and
\[
\limsup_{|x|\to\infty} \frac{j(-x)}{j(x)} < \infty,
\]
where $j:\erre\to\erre$ is a convex function such that $\beta=\partial
j$ (such a function always exists and is unique modulo addition of
constants -- see e.g. \cite{barbu-nonlin}, \cite{Bmax}). Here
$\partial$ stands for the subdifferential in the sense of convex
analysis.

\begin{definition}
  An adapted weakly c\`adl\`ag process $X \in L^1((0,T) \times \Xi \times
  \Omega)$ is a strong solution of (\ref{eq:spm}) if there exists an
  adapted process $\eta\in L^1((0,T) \times \Xi \times \Omega)$ such
  that $\eta(t,\xi) \in \beta(X(t,\xi))$ a.e. in $Q_T$, and
\begin{gather}
t \mapsto \int_0^t \eta(s,\xi)\,ds \in C^w([0,T],H_0^1),\\
X(t) - \Delta\int_0^t\eta(s)\,ds = x + \int_0^t B(X(s-))\,dM(s)
\quad \forall t\in[0,T],
\end{gather}
and $j(X)$, $j^*(\eta) \in L^1((0,T)\times\Xi\times\Omega)$. All
statements are meant to hold $\P$-a.s..
\end{definition}

We first establish a well-posedness result for the SPDE with additive
noise. Let us denote by $\mathcal{H}_2(T)$ and $\mathbb{H}_2(T)$ the
spaces of adapted processes $u:[0,T]\to H^{-1}$ such that
\[
\sup_{t\leq T} \E|u(t)|_{-1}^2 < \infty
\quad \text{and} \quad
\E\sup_{t\leq T} |u(t)|_{-1}^2 < \infty,
\]
respectively. We shall also denote by $\mathcal{H}_2$ the space of
$H^{-1}$-valued random variables with finite second moment.

The following intermediate result amounts to saying that the SPDE with
additive noise is globally well-posed. Note that we do not yet need to
assume that $M$ has stationary independent increments, as the latter
assumption will be used only in the proof of theorem \ref{thm:main}.

\begin{theorem}\label{thm:add}
  If $G\in\mathcal{G}(D((-\Delta)^\gamma))$, $\gamma>d/2$, then for
  each $x \in \mathcal{H}_2$ there exists a unique strong solution to
  the equation
\begin{equation}  \label{eq:bella}
\begin{cases}
dY(t)-\Delta\beta(Y(t))\,dt = G(t)\,dM(t) &\text{in $Q_T$}\\
\beta(Y(t))=0 &\text{in $\partial Q_T$}\\
Y(0)=x &\text{in $\Xi$}
\end{cases}
\end{equation}
Moreover, for $G_1$, $G_2\in \mathcal{G}(D((-\Delta)^\gamma))$ and
$y_1$, $y_2 \in \mathcal{H}_2$, denoting by $Y(t,y_i,G_i)$, $i=1,2$,
the solutions of (\ref{eq:bella}) with $G=G_i$ and $Y(0)=y_i$,
respectively, the following estimate holds:
$$
\E|Y(t,y_1,G_1)-Y(t,y_2,G_2)|_{-1}^2 \leq \E|y_1-y_2|_{-1}^2
+ \E\int_0^t |G_1(s)-G_2(s)|_{Q_M}^2\,d\langle M \rangle(s).
$$
Finally, the solution map $x \mapsto Y$ is a contraction from
$\mathcal{H}_2$ to $\mathbb{H}_2(T)$.
\end{theorem}

Our main result is the following.

\begin{theorem}
\label{thm:main}
Assume that $M$ has stationary independent increments. Then for each
$x\in \mathcal{H}_2$ there exists a unique strong solution of
(\ref{eq:main}). Moreover, the solution map $x \mapsto X$ is Lipschitz
from $\mathcal{H}_2$ to $\mathbb{H}_2(T)$.
\end{theorem}


\section{Auxiliary results}
Since the stochastic integral (\ref{eq:sint}) is a locally square
integrable martingale for any $G$ satisfying (\ref{eq:G}), Doob's
inequality yields the following simple result.
\begin{lemma}      \label{lem:GM}
  Let $M \in \mathcal{M}_{loc}^2(K)$ and $G \in \mathcal{G}(H)$.  Then
  \[
  \P\Big( \sup_{t\leq T} \big| G\cdot M(t) \big|_H < \infty \Big) = 1
  \]
\end{lemma}
\begin{proof}
It is enough to note that, by Cauchy-Schwartz' inequality,
\[
\E\sup_{t\leq T} \big| G\cdot M(t) \big|_H \leq
\big( \E\sup_{t\leq T} \big| G\cdot M(t) \big|^2_H \big)^{1/2},
\]
and that (since $|G\cdot M|_H^2$ is a submartingale) Doob's inequality
yields
\[
\E\sup_{t\leq T} \big| G\cdot M(t) \big|^2_H \leq
4 \E \big| G\cdot M(T) \big|^2_H =
4\E\int_0^T |G(s)|^2_{Q_M(s)}\,d\langle M\rangle(s)
       < \infty.
\]
The last step follows by the isometric formula (see \cite{MP-Z}) and
(\ref{eq:G}).
\end{proof}

\medskip

We shall also need an It\^o's formula for the square of the norm of
strong solutions to (\ref{eq:bella}).
\begin{lemma}\label{lm:itosq}
Let $Y$ be a strong solution of (\ref{eq:bella}).
Then one has
\begin{equation}
\label{eq:lazza}
\begin{split}
|Y(t)|_{-1}^2 &= |Y(0)|_{-1}^2 - 2\int_0^t \ip{Y(s)}{\eta(s)}_2\,ds\\
&\quad + 2\int_0^t \ip{Y(s-)}{G(s)\,dM(s)}_{-1} + [G \cdot M](s)\\
\end{split}
\end{equation}
for all $t\in[0,T]$, $\P$-a.s..
\end{lemma}
\begin{proof}
Let us set, for $m \in \mathbb{N}$ such that $m > 2 \vee (d+2)/4$,
\begin{align*}
  Y_\varepsilon(t) &= (1-\varepsilon \Delta)^{-m} Y(t),
  & \eta_\varepsilon &= (1-\varepsilon \Delta)^{-m} \eta(t),\\
  G_\varepsilon(t) &= (1-\varepsilon \Delta)^{-m} G(t),
  & x_\varepsilon &= (1-\varepsilon \Delta)^{-m} Y(0).
\end{align*}
Then we have by (\ref{eq:bella}) that
\[
dY_\varepsilon(t) = \Delta\eta_\varepsilon(t)\,dt + G_\varepsilon(t)\,dM(t),
\quad Y_\varepsilon(0)=x_\varepsilon, \quad
\eta_\varepsilon=0 \text{\ on } \partial\Xi,
\]
and It\^o's formula for $|Y_\varepsilon(t)|_{-1}^2$ yields
\begin{equation}
\label{eq:ito}
\begin{split}
  |Y_\varepsilon(t)|_{-1}^2 &= |x_\varepsilon|_{-1}^2
                - 2\int_0^t\ip{Y_\varepsilon(s-)}{\eta_\varepsilon(s)}_2\,ds
                + 2\int_0^t \ip{Y_\varepsilon(s-)}{G_\varepsilon(s)\,dM(s)}_{-1}\\
    &\quad + \int_0^t d\tr\tqv{G_\varepsilon\cdot M}^c(s)
           + \sum_{s \leq t} \big|\Delta(G_\varepsilon\cdot M)(s)\big|_{-1}^2\\
    &= |x_\varepsilon|_{-1}^2
       - 2\int_0^t\ip{Y_\varepsilon(s-)}{\eta_\varepsilon(s)}_2\,ds
       + 2\int_0^t \ip{Y_\varepsilon(s-)}{G_\varepsilon(s)\,dM(s)}_{-1}\\
    &\quad + [G_\varepsilon\cdot M](t)
\end{split}
\end{equation}
$\P$-a.s., where we have used the identity $\tr\tqv{Z}=[Z]$, which
holds for any semimartingale $Z$. We clearly have
$|x_\varepsilon|_{-1}^2 \uparrow |Y(0)|_{-1}^2$ $\P$-a.s. as $\varepsilon
\to 0$.  Moreover, we have
\begin{equation}    \label{eq:aprilia}
\int_0^t \ip{Y_\varepsilon(s-)}{\eta_\varepsilon(s)}_2\,ds =
\int_0^t \ip{Y_\varepsilon(s)}{\eta_\varepsilon(s)}_2\,ds \to
\int_0^t \ip{Y(s)}{\eta(s)}_2\,ds,
\end{equation}
as it follows from lemma 3.1 of \cite{BDPR-porous}. In fact, recalling
that $Y_\varepsilon(s)=y_\varepsilon(s)+G_\varepsilon\cdot M(s)$,
where $y_\varepsilon$ is weakly continuous and $G_\varepsilon\cdot M$
is c\`adl\`ag, we have
\[
\int_0^t \ip{y_\varepsilon(s-)}{\eta_\varepsilon(s)}_2\,ds =
\int_0^t \ip{y_\varepsilon(s)}{\eta_\varepsilon(s)}_2\,ds
\]
by weak continuity of $y_\varepsilon$, and
\[
\int_0^t \ip{G_\varepsilon\cdot M(s) - G_\varepsilon\cdot M(s-)}
{\eta_\varepsilon(s)}_2\,ds = 0
\]
because the times of discontinuity of c\`adl\`ag processes are at most
countable, hence with Lebesgue measure zero.

For the last term on the right hand side of
(\ref{eq:ito}) we can write
\[
\E[G_\varepsilon \cdot M](t) =
\E\langle G_\varepsilon \cdot M \rangle(t) \leq
\E\int_0^t |G_\varepsilon(s)|_{Q_M}^2\,d\langle M\rangle(s) \leq
\E\int_0^t |G(s)|_{Q_M}^2\,d\langle M\rangle(s) < \infty,
\]
hence by monotone convergence we have $\E[G_\varepsilon \cdot M](t)
\to \E[G\cdot M](t)$, and $[G_\varepsilon \cdot M](t) \to [G \cdot
M](t)$ $\P$-a.s. for all $t\in[0,T]$, passing to a subsequence if
necessary.

Let us now consider the third term on the right hand side of
(\ref{eq:ito}). We can write
\[
\begin{split}
&\E\Big|\int_0^t \ip{Y_\varepsilon(s-)}{G_\varepsilon(s)\,dM(s)}_{-1} -
\int_0^t \ip{Y(s-)}{G(s)\,dM(s)}_{-1} \Big|\\
&\quad \leq
\E\Big|\int_0^t \ip{Y(s-)}{(G_\varepsilon(s)-G(s))\,dM(s)}_{-1}
\Big|\\
&\qquad +
\E\Big|\int_0^t \ip{Y_\varepsilon(s-)-Y(s-)}{G_\varepsilon(s)\,dM(s)}_{-1}
\Big| =: I_1 + I_2.\\
\end{split}
\]
Davis' inequality and Cauchy-Schwartz' inequality yield (here we adopt
the usual notation $Y^*(t):=\sup_{s\leq t} |Y(s)|_{-1}$)
\begin{align*}
  I_1 &\leq \E\big[ Y_- \cdot((G_\varepsilon-G)\cdot M)\big]^{1/2}(t)
  \leq \E Y^*(t) [(G_\varepsilon-G) \cdot M]^{1/2}(t)\\
  &\leq \big(\E Y^*(t)^2\big)^{1/2}
        \big(\E [(G_\varepsilon-G) \cdot M](t)\big)^{1/2}\\
  &= \big(\E Y^*(t)^2\big)^{1/2}
        \big(\E \langle(G_\varepsilon-G) \cdot M\rangle(t)\big)^{1/2}\\
  &= \big(\E Y^*(t)^2\big)^{1/2}
     \Big( \E\int_0^t |G_\varepsilon(s)-G(s)|_{Q_M}^2 d\langle M\rangle(s)
     \Big)^{1/2},
\end{align*}
which converges to zero as $\varepsilon \to 0$ by (\ref{eq:G}) and
dominated convergence, provided we can show that $\E (Y^*(t))^2 <
\infty$. By (\ref{eq:ito}), recalling that
$(1-\varepsilon\Delta)^{-m}$ is a contraction, we have
\[
\E Y_\varepsilon^*(T)^2 \leq \E|Y(0)|^2 
+ 2\E\sup_{s\leq T} \Big|\int_0^s
      \ip{Y_\varepsilon(r-)}{G_\varepsilon(r)dM(r)}_{-1}\Big|
+ [G \cdot M](T),
\]
where we have used the inequality
\begin{equation}     \label{eq:spe}
\int_0^t \ip{Y_\varepsilon(s)}{\eta_\varepsilon(s)}_2\,ds \geq 0
\qquad \P\text{-a.s.}.
\end{equation}
The latter holds true by the following argument: since $\eta(t,x) \in
\beta(Y(t,x))$ for almost all $(t,x)\in Q_T$ $\P$-a.s., then
$Y(t,x)\eta(t,x) \geq 0$ a.e. in $Q_T$ $\P$-a.s. by monotonicity of
$\beta$. Since $(1-\varepsilon \Delta)^{-1}$ preserves the sign, thus
so does also $(1-\varepsilon\Delta)^{-m}$, one infers that
$Y_\varepsilon(t,x)\eta_\varepsilon(t,x) \geq 0$ a.e. in $Q_T$
$\P$-a.s., which implies (\ref{eq:spe}).
Setting $Z_\varepsilon=G_\varepsilon \cdot M$, Davis' inequality and
the elementary inequality $ab \leq \delta
a^2+\frac{b^2}{\delta}$ yield
\begin{align*}
&\E\sup_{t\leq T}
\left|\int_0^t \ip{Y_\varepsilon(s-)}{G_\varepsilon(s-)\,dM(s)}_{-1}\right|
= \E\sup_{t\leq T}
\left|\int_0^t \ip{Y_\varepsilon(s-)}{dZ_\varepsilon(s)}_{-1}\right|\\
 &\qquad\leq 3\E\big[Y_{\varepsilon-} \cdot Z_\varepsilon\big](T)^{1/2}
        \leq 3\E Y_\varepsilon^*(T) [Z_\varepsilon](T)^{1/2}\\
 &\qquad\leq 3\delta \E Y_\varepsilon^*(T)^2
       + \frac{3}{\delta} \E\int_0^T |G(s)|_{Q_M}^2\,d\langle M\rangle(s),
\end{align*}
because
\[
\E[Z_\varepsilon](T) = \E\langle Z_\varepsilon \rangle (T)
\leq \E \int_0^T |G(s)|_{Q_M}^2\,d\langle M \rangle(s).
\]
Thus we obtain
\begin{equation*}
\E\sup_{s\leq T}|Y_\varepsilon(s)|^2_{-1} \leq N
\E\int_0^T |G(s))|_{Q_M}^2\,d\langle M \rangle(s) < \infty,
\end{equation*}
where $N$ is a constant independent of $\varepsilon$.
Since $\varepsilon \mapsto \E Y_\varepsilon^*(T)^2$ is an increasing
bounded sequence, Fatou's lemma allows to conclude that $\E
Y^*(T)^2<\infty$, hence finally that $I_1 \to 0$ as $\varepsilon \to
0$.
A similar reasoning shows that $I_2 \to 0$ as $\varepsilon \to 0$.  We
thus have, by Chebishev's inequality, that
\[
\int_0^t \ip{Y_\varepsilon(s-)}{G_\varepsilon(s)\,dM(s)}_{-1} \to
\int_0^t \ip{Y(s-)}{G(s)\,dM(s)}_{-1}
\]
$\mathbb{P}$-a.s. and for all $t\in[0,T]$, at least on a subsequence of
$\varepsilon$, still denoted by $\varepsilon$, as $\varepsilon \to 0$.
\end{proof}


\section{Proof of theorem \ref{thm:add}}
The proof will be sketched only, underlying the differences with
respect to the corresponding proof in \cite{BDPR-porous}.

Let us consider the approximating SPDE (in integral form)
\begin{equation}
\label{eq:aspe}
X(t) = x + \int_0^t \Delta\big(\beta_\lambda(X(s))+\lambda X(s)\big)\,ds
+ G\cdot M(t),
\end{equation}
where $\beta_\lambda=\lambda^{-1}(I-(I+\lambda\beta)^{-1})$, $\lambda>0$, is the
Yosida approximation of $\beta$.

Then one has the following result.
\begin{lemma}
The SPDE (\ref{eq:aspe}) admits a unique c\`adl\`ag adapted solution
$X_\lambda$ such that
$$
X_\lambda, \; \beta_\lambda(X_\lambda) \in L^2([0,T],H_0^1).
$$
\end{lemma}
\begin{proof}
Equation (\ref{eq:aspe}) can be equivalently rewritten as the
deterministic PDE with random coefficients
\begin{equation}
\label{eq:dpde}
y'=\Delta\tilde\beta_\lambda(y+G\cdot M),  
\end{equation}
setting $y=X-G\cdot M$ and
$\tilde\beta_\lambda(x):=\beta_\lambda(x)+\lambda x$. Moreover, for
any fixed $\omega\in\Omega$, the time-dependent operator
\begin{align*}
  \mathcal{A}(t): H_0^1 &\to H^{-1}\\
   x &\mapsto -\Delta\tilde\beta_\lambda(x+G\cdot M)
\end{align*}
satisfies the assumptions of Theorem III.4.2 in \cite{barbu-nonlin},
hence (\ref{eq:dpde}) admits a unique solution
\[
y_\lambda\in C([0,T],L^2) \cap L^2([0,T],H_0^1),
\]
with $y'_\lambda \in L^2([0,T],H^{-1})$. Moreover, $y_\lambda$
depends continuously on $G\cdot M$ with respect to pathwise
convergence in $H^{-1}$, hence $X_\lambda:=y_\lambda+G\cdot M$ is an
adapted c\`adl\`ag solution of (\ref{eq:aspe}), as required.
\end{proof}
\begin{remark}
  Since $\beta_\lambda$ is Lipschitz, one can immediately conclude by
  \cite[Thm~24.7]{Met} that (\ref{eq:aspe}) has a unique c\`adl\`ag
  (strong) solution taking values in $H^{-1}$. It does not seem
  immediate to obtain also that the solution belongs to
  $L^2([0,T],H^1_0)$.
\end{remark}

We shall need some a priori estimates for
$z_\lambda:=(1+\lambda\beta)^{-1}X_\lambda$ and
$\eta_\lambda:=\beta_\lambda(X_\lambda)$.
\begin{lemma}
There exists $\Omega_0 \subset \Omega$ with $\P(\Omega_0)=1$ such
that, for every fixed $\omega \in \Omega_0$, one has
\begin{align}
\label{eq:313}
\int_{Q_T} (j(z_\lambda)+j^*(\eta_\lambda))\,d\xi\,ds &\leq N_1(1+|x|_{-1}^2),\\
\label{eq:314}
\int_{Q_T} |X_\lambda - z_\lambda|^2\,d\xi\,ds &\leq 2\lambda N_1 (1+|x|_{-1}^2),
\end{align}
where $N_1$ is a positive constant that may depend on $\omega$.
\end{lemma}
\begin{proof}
  By lemma \ref{lem:GM}, Sobolev's embedding theorem
  $D((-\Delta)^\gamma) \subset L^\infty$, $\gamma>d/2$, and the
  hypothesis that $\E|x|_{-1}^2<\infty$, it follows that there exists
  $\Omega_0 \subseteq \Omega$, $\P(\Omega_0)=1$, such that
  \[
  \sup_{t\leq T} |G\cdot M(t)|_{L^\infty} < \infty,
  \quad |x(\omega)|_{-1}^2 < \infty
  \qquad
  \forall \omega\in\Omega_0.
  \]
  Using this estimate in place of the corresponding one for $W_G$ in
  \cite[{\S}3.1]{BDPR-porous}, the claim follows.
\end{proof}
The above estimates allow us to pass to the limit as $\lambda \to 0$,
as in \cite[{\S}3.2]{BDPR-porous}, obtaining the following result.
\begin{lemma}
  There exist $y \in C^w([0,T],H^{-1}) \cap L^1(Q_T)$ and $\eta \in
  L^1(Q_T) \cap L^\infty([0,T],H_0^1)$ such that
\begin{equation}
  \label{eq:sati}
y(t) + A\int_0^t \eta(s)\,ds = x.
\end{equation}
\end{lemma}
One then continues proving that $\eta \in \beta(y+G\cdot M)$ a.e. in
$Q_T$, and that such $y$ and $\eta$ are unique. Hence the above
convergence results hold $\P$-a.s. for any choice of the sequence
$\lambda$. In particular, $y$ and $\eta$ are adapted processes.
Moreover, since $G\cdot M$ is c\`adl\`ag and $y$ is weakly continous,
it follows that $Y(t)=y(t)+G\cdot M(t)$ is an $H^{-1}$-valued weakly
c\`adl\`ag process such that
\[
Y(t) - \Delta \int_0^t \eta(s)\,ds = x + G \cdot M(t) 
\qquad \forall t \in [0,T]\quad \P\text{-a.s.},
\]
i.e. $Y$ solves (\ref{eq:bella}).

Once existence has been established, we need to prove uniqueness and
continuous dependence on the initial datum. This can be achieved with
the help of Lemma \ref{lm:itosq}. In particular, taking into account
that the second term on the right hand side of (\ref{eq:lazza}) is
negative because $\eta(s) \in \beta(Y(s))$ $\P$-a.s. for a.a. $s \in
[0,T]$, we have, by Lemma \ref{lm:itosq},
\begin{equation}      \label{eq:honda}
\begin{split}
|Y_1(t)-Y_2(t)|_{-1}^2 &\leq
            |y_1-y_2|_{-1}^2 
          + 2\int_0^t \ip{Y_1(s-)-Y_2(s-)}{(G_1(s)-G_2(s))\,dM(s)}_{-1}\\
   &\quad + [(G_1-G_2)\cdot M](t),
\end{split}
\end{equation}
where we set, for simplicity of notation, $Y_i:=Y(\cdot,y_i,G_i)$, $i=1,2$.
Taking expectation on both sides, we are left with
\begin{equation}     \label{eq:resta}
\E|Y_1(t)-Y_2(t)|_{-1}^2 \leq \E|y_1-y_2|_{-1}^2 
      + \E\int_0^t \big|(G_1(s)-G_2(s)\big|^2_{Q_M}\,d\langle M\rangle(s).
\end{equation}
Similarly, if $G_1=G_2$, (\ref{eq:honda}) immediately yields
\[
\E \sup_{t\leq T} |Y(t,y_1)-Y(t,y_2)|_{-1}^2 \leq
\E|y_1-y_2|_{-1}^2.
\]

\section{Proof of theorem \ref{thm:main}}
Consider the equation
\begin{equation}
\label{eq:tos}
dY(t)=\Delta\beta(Y(t))\,dt+B(X(t-))\,dM(t),
\quad t\in[0,T],
\end{equation}
and define the operator $\Phi:X \mapsto Y$ that associates to
$X\in\mathbb{H}_2(T)$ the solution $Y$ of (\ref{eq:tos}). We are going
to prove that $\Phi$ is an endomorphism of $\mathbb{H}_2(T)$ and is a
contraction. Moreover, since $t\mapsto B(X(t-))$ is predictable, we
know by theorem \ref{thm:add} that $Y$ is adapted and weakly
c\`adl\`ag.  Let us first obtain two estimates that hold for any
quasi-left-continuous $M \in \mathcal{M}^2_{loc}(K)$. It\^o's formula
yields
\begin{multline*}
|Y(t)|_{-1}^2 + 2\int_0^t\ip{Y(s)}{\eta(s)}_2\,ds =\\
|Y(0)|_{-1}^2 + 2\int_0^t \ip{Y(s-)}{B(X(s-))\,dM(s)}_{-1}
+ \big[B(X_-)\cdot M](t),
\end{multline*}
where $B(X_-)$ stands for $t\mapsto B(X(t-))$.
Since $\ip{Y(s)}{\eta(s)} \geq 0$ $\P$-a.s. for all $s\leq t$, we can write
\begin{equation}
\label{eq:stima0}
\E\sup_{t\leq T} |Y(t)|_{-1}^2 \leq 
2 \E\sup_{t\leq T} \left|\int_0^t \ip{Y(s-)}{(B(X(s-))\,dM(s)}_{-1}\right|
+ \E \int_0^T |B(X(s)|_{Q_M}^2\,d\langle M \rangle(s).
\end{equation}
Following in a completely similar way as we have done in the last
paragraph of the proof of Lemma \ref{lm:itosq}, we obtain
\begin{equation}
\label{eq:cano}
(1-6\varepsilon) |Y|^2_{\mathbb{H}_2(T)} \leq
(6/\varepsilon + 1) \E\int_0^T |B(X(s))|_{Q_M}^2\,d\langle M \rangle(s).
\end{equation}
Similarly, writing
\[
\left\{
\begin{array}{l}
\ds dY_1 = \Delta\beta(Y_1)\,dt + B(X_{1-})\,dM\\[4pt]
\ds dY_2 = \Delta\beta(Y_2)\,dt + B(X_{2-})\,dM,
\end{array}
\right.
\]
with $Y_1(0)=Y_2(0)$, using again It\^o's formula (see lemma
\ref{lm:itosq}), in complete analogy to the above derivation, we
obtain the estimate
\begin{equation}
\label{eq:gari}
(1-6\varepsilon) |Y_1-Y_2|^2_{\mathbb{H}_2(T)} \leq
(6/\varepsilon + 1) 
\E\int_0^T |B(X_1(t))-B(X_2(t))|_{Q_M}^2\,d\langle M \rangle(t).
\end{equation}

If $M \in \mathcal{M}^2_{loc}(K)$ has also stationary independent
increments, then, in view of remark \ref{rmk:pissi},
we have
\begin{eqnarray*}
\E\int_0^T |B(X(s))|_{Q_M}^2\,d\langle M\rangle(s) &=&
\E\int_0^T |B(X(s))|_Q^2\,ds \\
&\leq& k\E\int_0^T(1+|X(s)|^2)\,ds\\
&\leq& kT \big(1+|X|^2_{\mathbb{H}_2(T)}\big) < \infty,
\end{eqnarray*}
hence, choosing $\varepsilon < 1/6$, $|Y|^2_{\mathbb{H}_2(T)}<\infty$,
by virtue of (\ref{eq:cano}). This proves that the image of $\Phi$ is
contained in $\mathbb{H}_2(T)$. Let us now show that $\Phi$ is a
contraction. In fact, (\ref{eq:gari}) and assumption (\ref{eq:assb})
yield
$$
|Y_1-Y_2|^2_{\mathbb{H}_2(T)} \leq
\frac{1+6/\varepsilon}{1-6\varepsilon}\,kT
|X_1-X_2|^2_{\mathbb{H}_2(T)},
$$
i.e. $\Phi$ is a contraction on $\mathbb{H}_2(T)$ whenever
\begin{equation}
  \label{eq:smallt}
T < \frac{1-6\varepsilon}{1+6/\varepsilon} \, \frac1k.
\end{equation}
Then, by the Banach fixed point theorem, there exists a unique solution of
(\ref{eq:spm}). If $T$ does not satisfy (\ref{eq:smallt}), then one
proceeds in a classical way considering intervals $[0,T_0]$,
$[T_0,2T_0]$, etc., with suitably small $T_0$, such that $\Phi$ is a
contraction on $\mathbb{H}_2(T_0)$.

In order to prove Lipschitz continuity of the solution map, note that
we have
\[
(1-6\varepsilon)|Y(\cdot,y_1)-Y(\cdot,y_2)|^2_{\mathbb{H}_2(T)}
\leq
(6/\varepsilon+1) k T |Y(\cdot,y_1)-Y(\cdot,y_2)|^2_{\mathbb{H}_2(T)}
+|y_1-y_2|^2_{\mathcal{H}_2},
\]
hence for any $T_0$ such that
\[
1 - 6\varepsilon - k T_0 (6/\varepsilon + 1) > 0
\]
we have
\[
|Y(\cdot,y_1)-Y(\cdot,y_2)|_{\mathbb{H}_2(T_0)} \leq
N_0 |y_1-y_2|_{\mathcal{H}_2},
\]
where
\[
N_0 = \big(
1 - 6\varepsilon - k T_0 (6/\varepsilon + 1)
\big)^{-1/2}.
\]
Considering intervals of length $T_0$ covering $[0,T]$ one finally gets
\[
|Y(\cdot,y_1)-Y(\cdot,y_2)|_{\mathbb{H}_2(T)} \leq
N |y_1-y_2|_{\mathcal{H}_2},
\]
where $N=N(k,T)$.


\section{Generalized solutions}
In this section we introduce a concept of generalized solution for
equation (\ref{eq:spm}), which allows to replace the assumption
(\ref{eq:asso}) with
\begin{equation}     \label{eq:ace}
B: H^{-1} \to \mathcal{L}_2^Q(K,H^{-1}).
\end{equation}
As we did before, we start with the case of additive noise and general
$M \in \mathcal{M}_{loc}^2(K)$.
\begin{definition}
  Let $G \in \mathcal{G}(H^{-1})$. An adapted process $Y$ is called a
  $\mathcal{H}$-generalized solution of (\ref{eq:bella}) if there
  exists a sequence $\{G_n\}_{n \in \mathbb{N}} \subset
  \mathcal{G}(D((-\Delta)^\gamma))$ with
  \[
  \lim_{n\to\infty}
     \E \int_0^T |G_n(s)-G(s)|^2_{Q_M}\,d\langle M \rangle(s) = 0
  \]
  such that the solution $Y_n$ to
  \[
  dY(t) = \Delta\beta(Y(t))\,dt + G_n(t)\,dM(t),
  \]
  equipped with the same initial and boundary conditions of
  (\ref{eq:bella}), converges to $Y$ in $\mathcal{H}_2(T)$. If the
  convergence is in $\mathbb{H}_2(T)$, $X$ is called
  $\mathbb{H}$-generalized solution.
\end{definition}
It is clear that a $\mathbb{H}$-generalized solution is also a
$\mathcal{H}$-generalized solution. In the following we shall refer to
$\mathbb{H}$-generalized solutions simply as generalized solutions.

\begin{theorem}     \label{thm:gener}
  Let $G \in \mathcal{G}(H^{-1})$. Then (\ref{eq:bella}) admits a
  unique generalized solution. Moreover, the solution map $x \mapsto
  Y$ is a contraction from $\mathcal{H}_2$ to $\mathbb{H}_2(T)$.
\end{theorem}
For the proof of the theorem we need the following approximation
procedure for elements of the space $H^{-1}$.  Let $f \in
H^{-1}$. Then there exists $F \in H^1_0$ such that $f=\Delta F$, and
$|f|_{-1}^2 = |\nabla F|_2^2$. Set $F_n = \zeta_n \ast F$ and $f_n =
\Delta F_n$, where $\{\zeta_n\}_{n \in \mathbb{N}}$ is a standard
sequence of mollifiers (here we have considered an extension of $F$ to
$H^{-1}(\erre^d)$, still denoted by $F$). In particular, since $F_n
\in C^\infty$, then $f_n \in C^\infty(\Xi) \subset L^\infty(\Xi)$.
Recalling that $\nabla(\zeta_n\ast F) = \zeta_n\ast\nabla F$, we have
\begin{equation}     \label{eq:amatri}
|f_n-f|^2_{-1} = |\nabla F_n - \nabla F|^2_2 =
|\zeta_n \ast (\nabla F) - \nabla F|_2^2 \to 0
\end{equation}
as $n \to \infty$ because $\nabla F \in L^2$ and $\zeta_n \ast \phi
\to \phi$ in $L^2$ for all $\phi \in L^2$.
Moreover, Young's inequality for convolutions yields
\begin{equation}    \label{eq:ciana}
|f_n|_{-1} = |\zeta_n \ast \nabla F|_2
\leq |\zeta_n|_1 |\nabla F|_2 \leq |\nabla F|_2 = |f|_{-1}.
\end{equation}
The map associating $f$ to $f_n$ will de noted by $\Lambda_n$.

\begin{proof}[Proof of theorem \ref{thm:gener}]
Uniqueness follows by (\ref{eq:resta}). In fact, the estimate is
stable with respect to passage to the limit in $G_1$ and $G_2$.

For fixed $s\in [0,T]$ and $\kappa \in Q_M^{1/2}K$ we have $G(s)\kappa
\in H^{-1}$, and we define $G_n(s)\kappa := \Lambda_n G(s)\kappa$.
Let us show that
\[
\lim_{n \to \infty} \E\int_0^T |G_n(s)-G(s)|^2_{Q_M}\,
d\langle M\rangle(s) = 0.
\]
The claim follows by (\ref{eq:amatri}), the dominated convergence
theorem and the inequality
\[
|G_n(s)Q_M^{1/2}e_k|_{-1} \leq |G(s)Q_M^{1/2}e_k|_{-1}
\qquad \P\text{-a.s.},
\]
which holds for all $s\in[0,T]$ and all $k\in\mathbb{N}$, where
$(e_k)_{k\in\mathbb{N}}$ is a basis of $K$, as it follows by
(\ref{eq:ciana}).

\noindent
Now (\ref{eq:resta}) implies that
\[
\sup_{t\leq T} \E |Y_n(t)-Y_m(t)|_{-1}^2 \leq
\E\int_0^T |G_n(s)-G_m(s)|^2_{Q_M}\,
d\langle M\rangle(s),
\]
that is $\{Y_n\}_{n\in\mathbb{N}}$ is a Cauchy sequence in
$\mathcal{H}_2(T)$, which converges to a $\mathcal{H}$-generalized
solution $Y$ of (\ref{eq:bella}).

Let us show that $\{Y_n\}_{n\in\mathbb{N}}$ is a Cauchy sequence also
in $\mathbb{H}_2(T)$, which proves the existence of a
$\mathbb{H}$-generalized solution. In fact, using an argument based on
It\^o's formula and Davis' inequality completely analogous to tho one
leading to (\ref{eq:gari}), we obtain
\[
\E \sup_{t\leq T} |Y_n(t)-Y_m(t)|^2 \leq
N \E \int_0^T |G_n(t)-G_m(t)|_{Q_M}^2\,d\langle M \rangle(t),
\]
where $N$ is a positive constant.
\end{proof}

It is now possible to extend the result to equations with
multiplicative noise.
\begin{theorem}
  Assume that $M$ has stationary independent increments and $B$ is as
  in (\ref{eq:ace}). Then (\ref{eq:main}) admits a unique generalized
  solution. Moreover, the solution map $x \mapsto X$ is Lipschitz
  from $\mathcal{H}_2$ to $\mathbb{H}_2(T)$.
\end{theorem}
\begin{proof}
  The argument is an extension of that used in the proof of theorem
  \ref{thm:main}, using the previous theorem. In fact, let $X \in
  \mathbb{H}_2(T)$ and consider equation (\ref{eq:tos}), which admits
  a unique generalized solution by theorem \ref{thm:gener}. Since
  estimates (\ref{eq:cano}) and (\ref{eq:gari}) hold also for
  generalized solutions (by a now obvious limiting procedure), the map
  associating $Y$ to $X$, as defined in the proof of theorem
  \ref{thm:main}, is a contraction in $\mathbb{H}_2(T_0)$ for a
  suitably small $T_0$. The rest of the proof is identical to that of
  theorem \ref{thm:main}.
\end{proof}

\bibliographystyle{amsplain}
\bibliography{ref}

\end{document}